\documentclass[12pt]{amsart}
\usepackage{amsmath,amssymb,amsthm,amscd}

\newtheorem{theorem}{Theorem}[section]

\theoremstyle{definition}

\theoremstyle{remark}

\begin{document}
 
\title[Rigidity of entire self-shrinking solutions]
{Rigidity of Entire self-shrinking solutions to curvature flows}

\author{Albert CHAU}
\address{Department of Mathematics\\
University of British Columbia\\
Vancouver, B.C., V6T 1Z2\\
Canada}
\address{Department of Mathematics \\University of Washington\\Seattle, WA 98195\\U.S.A.}
\email{chau@math.ubc.ca}

\author{Jingyi CHEN}
\email{jychen@math.ubc.ca}

\author{Yu Yuan }
\email{yuan@math.washington.edu}
\thanks{The first two authors are partially supported by NSERC, and the third author is
partially supported by NSF}
\date{\today}

\begin{abstract}
We show that (a) any entire graphic self-shrinking solution to the Lagrangian mean 
curvature flow in ${\mathbb C}^{m}$ with the Euclidean metric is flat; (b) any space-like entire graphic 
self-shrinking solution to the Lagrangian mean curvature flow in ${\mathbb C}^{m}$ with the pseudo-Euclidean metric 
is flat if the Hessian of the potential is bounded below quadratically; and (c) the Hermitian counterpart of 
(b) for the K\"ahler Ricci flow.
\end{abstract}

\maketitle

\section{Introduction}

Self-similar solutions to  curvature flows play an important role in understanding the general behavior of the flow and the types of singularities that can develop. For mean curvature flow, self-shrinking 
solutions arise naturally at a so-called type-I singularity from Huisken's monotonicity formula \cite{Huisken}. More precisely,  these are ancient families of immersions ${F}(x, t):\Sigma\times(-\infty, 0) \to{\mathbb R}^N$ of some manifold $\Sigma$ into ${\mathbb R}^N$ which solve the mean curvature flow equation 
\begin{equation}\label{MCF}
\left(\dfrac{d}{dt} { F }\right)^\perp={ H}
\end{equation}
simply by scaling  ${ F}(\Sigma, t)=\sqrt{-t}\,{ F}(\Sigma, -1)$.  Here $\left(\frac{d}{dt} { F }\right)^\perp$ is the normal component of the vector $\frac{d}{dt} { F }$ and $H$ is the mean curvature of $F(\Sigma,t)$.  It follows that ${ F}(x, -1)$ satisfies equation
\begin{equation}\label{MCFshrinker}
{ H} + \frac{1}{2}\, { F } ^\perp=0.
\end{equation}
 Conversely, if an embedding ${ F}$ satisfies \eqref{MCFshrinker} then the corresponding solution to the mean curvature flow will be a self-shrinking solution.  If a Lagrangian graph 
$\{(x,Du(x)): x\in{\mathbb R}^n\}$ in ${\mathbb R}^{2n}$ satisfies (\ref{MCFshrinker}), then up to an additive constant the potential function $u$ solves 
\begin{equation}\label{LMCF}
\arctan \lambda_1(x) + \cdot\cdot\cdot + \arctan \lambda_n (x) = \frac{1}{2} x\cdot Du(x)-u(x)
\end{equation}
where $\lambda_1(x), ..., \lambda_n(x)$ are the eigenvalues of the Hessian $D^2 u$ of $u$ at $x\in{\mathbb R}^n$.  The first main result in this note is the following Bernstein type rigidity for entire self-shrinking solutions for Lagrangian mean curvature flow: 
\begin{theorem}\label{maintheorem1}
If  $u(x)$ is an entire smooth solution to equation (\ref{LMCF}) in ${\mathbb R}^n$, then $u(x)$ is the quadratic polynomial $u(0)+\frac{1}{2}\langle D^2u(0)\,x,x\rangle$. 
\end{theorem} 

When ${\mathbb R}^{2n}$ is equipped with the indefinite metric $\sum_idx_i dy_i$,  if a space-like gradient graph $\{(x,Du(x)): x\in{\mathbb R}^n\}$ satisfies  (\ref{MCFshrinker}) then up to an additive constant the potential $u$ is convex and satisfies the elliptic equation

\begin{equation}\label{non-euclid}
\ln\det D^2 u(x) =\frac{1}{2} x\cdot Du(x) -u(x). 
\end{equation}
We have

\begin{theorem}\label{thm2}
If $u(x)$ is an entire smooth convex solution to (\ref{non-euclid}) in $\mathbb{R}^{n}$, then $u(x)$ is the quadratic polynomial $u(0)+\frac{1}{2}\langle D^2u(0)\,x,x\rangle$, provided either 
\begin{enumerate}
\item [(i)]  $\displaystyle{D^2 u (x)\geq \frac{2(n-1+\delta)}{|x|^2}}$ for any  $\delta>0$ as $|x|\to\infty$ or 
 \item [(ii)] $u$ is radially symmetric. 
\end{enumerate}
\end{theorem}

Finally, we also consider the Hermitian analog of (\ref{non-euclid}) and Theorem \ref{thm2}.   Namely we consider real valued functions 
satisfying  
\begin{equation}\label{eqn:10}
\ln\det \partial\bar{\partial} v(x) =\frac{1}{2}x\cdot D v(x) -v(x)
\end{equation}
on ${\mathbb C}^m$. This  is closely related to the K\"ahler Ricci flow as we describe in Section 4.  We prove
\begin{theorem}\label{thm3}
If $v(x)$ is an entire smooth pluri-subharmonic solution to (\ref{eqn:10}) in $\mathbb{C}^m$, then $v(x)$ is the quadratic polynomial $v(0)+\frac{1}{2}\langle D^2v(0)\,x,x\rangle_{\mathbb{R}^{2m}}$, provided either
\begin{enumerate}
\item [(i)]  $\displaystyle{\partial\bar{\partial}v(x)\geq \frac{2m-1+\delta}{2|x|^2} \, I}$  for any $\delta>0$ as $|x|\to\infty$ or 
 \item [(ii)] $v$ is radially symmetric. 
\end{enumerate}
\end{theorem}

By using existence \cite{CCH1} and uniqueness \cite{CP} results for the Lagrangian mean curvature flow, 
 the rigidity of self-expanding, self-shrinking and translating solutions for the Lagrangian mean curvature flow was studied in \cite{CCH2} when 
 the Hessian of the potential function is strictly bounded between $-1$ and $1$. The same 
rigidity for self-shrinking and translating solutions with arbitrarily bounded Hessian 
was derived from a Liouville type property for ancient solutions to parabolic equations 
\cite{NY} (for self-shrinking solutions, a special case of \cite{NY} was treated recently in \cite{HW2}). Theorem 
\ref{maintheorem1} improves the previous results on self-shrinking solutions by dropping  the assumption on 
Hessian completely. For the pseudo-Euclidean case, under a similar assumption on the Hessian as in Theorem 
\ref{thm2}, namely quadratically decaying lower bound, the Bernstein type result was obtained in \cite{HW}. 
However, our method is completely different and much simpler, and after scaling to the same equation (\ref{non-euclid}) it also gives a little sharper constant in the 
assumption on the Hessian.  Recently, it was shown in \cite{LW}  that any entire graphic hypersurface solution to \eqref{MCFshrinker} must be flat, hence generalizing an earlier result in \cite{EH}. 

A key ingredient in our arguments, for each of the three cases above, is that a natural geometric quantity, involving second order derivatives of the potential function, obeys a second order elliptic equation with an ``amplifying force".  We then construct a barrier function to show that the quantity is constant via the maximum principles. The homogeneity form of the lower order terms in the equations implies the potentials are quadratic polynomials. 

\section{Proof of Theorem \ref{maintheorem1}}
 We point out that if $u$ satisfies (\ref{LMCF}) then $v(x,t) = -t u(\frac{x}{\sqrt{-t}})$ satisfies
\begin{equation}\label{LMCF2}
\frac{\partial v}{\partial t} = -\sqrt{-1}\log\frac{\det\left(I+\sqrt{-1}D^2v\right)}{\sqrt{\det\left(I+D^2v D^2 v\right)}} 
\end{equation}
on ${\mathbb R}^{n}\times(-\infty, 0)$ and the family of embeddings ${F}$$(x, t) = (x,Dv(x,t))$ from ${\mathbb R}^{n}$ into ${\mathbb R}^{2n}$ solves the mean curvature flow \eqref{MCF} (cf. \cite{Sm}).  While this connection is our main motivation to study \eqref{LMCF}, we will not use \eqref{LMCF2} explicitly in our following proof of Theorem \ref{maintheorem1}.

Let $z_j = x_j +\sqrt{-1}y_j$ be the standard complex coordinates on ${\mathbb C}^n={\mathbb R}^{2n}$. The phase function 
$\Theta$ on a Lagrangian submanifold $\Sigma^n$, is defined by 
$$
dz_1\wedge ... \wedge dz_n |_{\Sigma^n} = e^{\sqrt{-1}\Theta}d\mu_{\Sigma^n}. 
$$ 
When $\Sigma^n$ is a Lagrangian graph $\{(x,Du(x)):x\in{\mathbb R}^n\}$ in ${\mathbb C}^n$,  $\Theta$ takes the form 
\begin{equation}\label{theta}
\Theta=\arctan\lambda_1+ ... +\arctan\lambda_n.
\end{equation}
For simplicity, for a function $f$, we denote $\frac{\partial^k f}{\partial x_{i_1} ... \partial x_{i_l}}$ by $f_{i_1 ... i_l}$ for $k=i_1+...+i_l$. Let 
\begin{eqnarray*}
A&=& \left(A_{ij}\right)=I+\sqrt{-1}D^2u,\,\,\,\left(A^{ij}\right)=A^{-1},\\
B&=&\left(B_{ij}\right)=I+D^2u\,D^2u,\,\,\,\left(B^{ij}\right)=B^{-1}.
\end{eqnarray*}
Observe that 
\begin{equation}
B =\left(I - \sqrt{-1}D^2u\right) A\,\,\,\,\,\mbox{and}\,\,\,\,\,\,A^{-1} = B^{-1}\left( I -\sqrt{-1}D^2 u\right)
\end{equation}
and 
\begin{equation}\label{Langle}
\Theta=-\sqrt{-1}\log \frac{\det A}{\sqrt{\det B}}. 
\end{equation} 
Then by \eqref{Langle} we have 
\begin{eqnarray}\label{1}
 \Theta_k &=& -\sqrt{-1}\sum_{i,j}\left( A^{ij} A_{ij,k} -\frac{1}{2}B^{ij}B_{ij,k}\right)\nonumber\\
 &=& -\sqrt{-1}\sum_{i,j,l}\left( B^{il}(\delta_{lj}-\sqrt{-1}u_{lj}) \cdot\sqrt{-1} u_{ijk}-\frac{1}{2}B^{ij} \cdot 2 u_{il}u_{ljk} \right)\nonumber\\
 &=&\sum_{i,j}B^{ij}u_{ijk}-\sqrt{-1} \sum_{i,j,l}\left( B^{il} u_{lj}u_{ijk}-B^{ij}u_{il}u_{ljk} \right) \nonumber\\
 &=&\sum_{i,j}B^{ij}u_{ijk}
 \end{eqnarray}
 as $B$ is symmetric and by changing indices.

Differentiating equation (\ref{LMCF}), and using \eqref{theta} we have
\begin{equation}\label{-1}
\Theta_i = -\frac{1}{2}u_i + \frac{1}{2} x\cdot D u_i
\end{equation}
and
\begin{equation}\label{0}
\Theta_{ij}= \frac{1}{2} x\cdot Du_{ij} =\frac{1}{2}\sum_k x_k u_{ijk}.
\end{equation}
Note that $B$ is just the induced metric $g$ of $\Sigma^n$ in ${\mathbb C}^n$ with the Euclidean metric. It follows from (\ref{0}) and (\ref{1}) that $\Theta$ satisfies   the following elliptic equation of non-diveregnce form: 
\begin{equation}\label{eqn:angle}
\sum_{i,j}g^{ij}\Theta_{ij}(x) - \frac{1}{2} x\cdot D\Theta(x) = 0
\end{equation}
with the ``amplifying force" $\frac{1}{2}x\cdot D\Theta(x)$. 

Next, we construct a radially symmetric barrier to show $\Theta$ attains its global maximum at an interior point. Take a radially symmetric function
\begin{equation}\label{solution}
w(r)=\epsilon\, r^{1+\delta}+\max_{\partial B_{r_0}}\{\Theta\}
\end{equation}
where $\epsilon$ is a positive constant and $B_{r_0}$ is the ball in ${\mathbb R}^n$ centered at the origin with radius $r_0=\sqrt{2(n-1+\delta)}$. 
For $|x|=r\geq r_0$, we have  
\begin{equation}\label{ode}
w_{rr}+ \frac{n-1}{r}w_r -\frac{r}{2}w_r\leq 0
\end{equation}  
where 
\begin{eqnarray*}
w_r &=&\epsilon \,(1+\delta)\, r^{\delta} > 0, \\
w_{rr}&=&\frac{\delta}{r}\,w_r>0.
\end{eqnarray*}
Also note
\[
D^{2}w\sim\left(
\begin{array}[c]{cccc}
w_{rr} &  &  & \\
& \frac{w_{r}}{r} &  & \\
&  & \cdots & \\
&  &  & \frac{w_{r}}{r}
\end{array}
\right)  
\geq 0 
\]
when $r> 0$. Observe that
\begin{equation}
g^{-1} = \left( I +D^2 u D^2 u\right)^{-1} \leq I.
\end{equation} 
Thus we have 
\begin{equation}
\mbox{tr}  \left(g^{-1} D^2 w \right) \leq \mbox{tr} \left( I D^2 w \right).
\end{equation} 
Hence, when $|x|\geq r_0$ we have 
\begin{equation}
\sum_{i,j}g^{ij} w_{ij} -\frac{1}{2} x\cdot Dw \leq \Delta w -\frac{1}{2} x\cdot Dw \leq  0
\end{equation}
where $ \Delta $ is the Euclidean Laplacian on ${\mathbb R}^n$ and we have used \eqref{ode} in the last equality.  

So far we have
\begin{equation}
\sum_{i,j}g^{ij}w_{ij}-\frac{1}{2}x\cdot Dw \leq \sum_{i,j}g^{ij}\Theta_{ij}-\frac{1}{2} x\cdot D\Theta\,\,\,\,\,\,\,\mbox{if $r_0
\leq |x| <\infty$} 
\end{equation}
with comparison along the boundaries:
\begin{equation}
w=\epsilon\, r_0^{1+\delta}+\max_{ \partial B_{r_0}}\{\Theta\}\geq \Theta \,\,\,\,\,\,\mbox{on $\partial B_{r_0}$ }
\end{equation} 
and
\begin{equation}
w(|x|)>\Theta (x)\,\,\,\,\,\,\,\,\,\,\,\,\,\,\, \mbox{when $|x|\to\infty,$}
\end{equation} 
since $\Theta$ is bounded  while $w(|x|)\to\infty$ as $|x|\to\infty$.  By the weak maximum principle, we get
\begin{equation}
\epsilon\, |x|^{1+\delta} +\max_{\partial B_{r_0}}\{\Theta\} = w(|x|) \geq \Theta(x) 
\end{equation}
for all $|x|\geq r_0$.  By letting $\epsilon$ go to zero, we then conclude that $\Theta$ achieves its global maximum on ${\mathbb R}^n$ in the closure of the ball $B_{r_0}$.  Applying the strong maximum principle to (\ref{eqn:angle}), we immediately see that $\Theta$ is a constant.  

Now from (\ref{0}), for any $i,j$ we have 
\begin{equation}
x\cdot D u_{ij}=0. 
\end{equation}
Euler's homogeneous function theorem asserts that $u_{ij}$ is homogenous of degree 0. However, the function $u_{ij}$ is smooth at the origin, therefore $u_{ij}$ is constant. It follows from (\ref{-1}) that $u$ is the quadratic polynomial in the claimed form.  This completes the proof of Theorem \ref{maintheorem1}.

\vspace{.3cm}

\noindent{\it Remark.} Denote $\widetilde{F}(x,t) = (x,Dv(x,t))$. Then by recalling $v(x, t)=\sqrt{-t} u (\frac{x}{ \sqrt{-t}})$ we have that (\ref{eqn:angle}) is equivalent to 
\begin{equation}\label{theta:parabolic}
\frac{\partial \Theta(\widetilde{F})}{\partial t} = \sum_{i,j}g^{ij}(\widetilde{F})\Theta(\widetilde{F})_{ij}.
\end{equation}
It is known that if $F(x,t)$ satisfies the mean curvature flow equation $F_t = H(F)$ and $F$ is Lagrangian, then its phase function $\Theta$ satisfies 
$$
\frac{\partial\Theta(F)}{\partial t}=\Delta_g \Theta(F)
$$ 
where $\Delta_g$ is the Laplace operator of the induced metric $g$ on the time slice $F(\cdot,t)$. The non-divergence structure of (\ref{theta:parabolic}) is due to the fact that 
$\widetilde{F}$ satisfies the mean curvature flow equation up to tangential diffeomorphisms.

\section{proof of theorem \ref{thm2}} 
We note that if $u$ satisfies (\ref{non-euclid}) then $v(x,t) = -t u(\frac{x}{\sqrt{-t}})$ verifies 
\begin{equation}\label{non-euclid2}
\frac{\partial v}{\partial t} =\ln\det D^2 u(x)
\end{equation}
on ${\mathbb R}^{n}\times(-\infty, 0)$ and the family of embeddings ${F}$$(x, t) = (x,Dv(x,t))$ from 
${\mathbb R}^{n}$ into ${\mathbb R}^{2n}$ solves the mean curvature flow \eqref{MCF} with respect to the  
pseudo-Euclidean background metric $ds^2 = \sum_idx^i dy^i$ on ${\mathbb R}^{2n}$ (cf. \cite{HW}).  Again, while 
this connection is our main motivation to study \eqref{non-euclid}, we will not use \eqref{non-euclid2} explicitly in our following proof of Theorem \ref{thm2}.

From (\ref{non-euclid}), we see that $D^2u>0$. Set $\Psi = \ln \det D^2 u$.  We have 
\begin{eqnarray*}
\Psi_i =\sum_{k,l}u^{kl}u_{kli} = \sum_{k,l}g^{kl}u_{kli}
\end{eqnarray*}
where $g^{-1}$ is the inverse of the induced metric $g$ of the graph $(x,Du(x))$ in ${\mathbb R}^{2n}$ with the pseudo-Euclidean metric above. 
On the other hand, by differentiating equation (\ref{non-euclid}) twice we obtain
\begin{eqnarray*}
\Psi_{ij}(x) = \frac{1}{2} x\cdot Du_{ij} (x)
\end{eqnarray*}
and hence, 
\begin{equation}\label{psi}
\sum_{i,j}g^{ij} \Psi_{ij}(x) -\frac{1}{2} x\cdot D\Psi (x)=0. 
\end{equation}

Next, as in the previous section, for any $\epsilon > 0$ we take a radially symmetric function $w$ defined by 
\begin{equation}
w(r) = \epsilon\, r^{1+\delta} + \max_{\partial B_1} \{-\Psi\}.
\end{equation}
It is clear that for $r$ positive 
\begin{eqnarray*}
w_r &=& \epsilon \,(1+\delta)\,r^\delta>0, \\
w_{rr}&=& \frac{\delta}{r} \,w_r >0,  
\end{eqnarray*}
\begin{equation}\label{eqn:w}
\frac{r^2}{2(n-1+\delta)}\,\Delta w-\frac{r}{2}w_r=0,
\end{equation}
where $\Delta$ is the Euclidean Laplacian on ${\mathbb R}^n$ and 
$$
D^2 w>0.
$$
By assumption (i) in Theorem \ref{thm2}, we have 
\begin{equation}
g^{-1} = (D^2 u )^{-1} \leq \frac{|x|^2}{2(n-1)+2\delta} \,I.
\end{equation}
Here we assume $\displaystyle{ D^2u\geq \frac{2(n-1+\delta)}{|x|^2}}$ for $|x|>1$ instead of $|x|$ being greater than a large number as in the assumption for simplicity. 
Otherwise, we just replace 1 by the large number, and our arguments go through as well.
Thus
\begin{equation}
\mbox{tr} \left(g^{-1} D^2 w\right) \leq  \frac{r^2}{2(n-1)+2\delta}\,\mbox{tr} \left( I D^2 w \right) 
\end{equation} 
and it follows that
\begin{equation}
\sum_{i,j}g^{ij}w_{ij} -\frac{1}{2} x \cdot Dw \leq \frac{r^2}{2(n-1)+2\delta}\,\Delta w   -\frac{1}{2}\, r \,w_r = 0
\end{equation}
where we have used (\ref{eqn:w}) to conclude the last equality. 

Thus far, we have 
\begin{equation}
\sum_{i,j}g^{ij}w_{ij} -\frac{1}{2} x \cdot Dw \leq \sum_{i,j} g^{ij}(-\Psi)_{ij}-\frac{1}{2}x\cdot D(-\Psi).
\end{equation} 
Also, we have that along the boundaries
$$
w(|x|)=\epsilon+\max_{\partial B_1}\{-\Psi\}\geq -\Psi(x)\,\,\,\,\,\,\,\mbox{on $\partial B_1$}
$$
and 
$$
w(|x|)> -\Psi(x)\,\,\,\,\,\,\,\,\mbox{as $|x|\to\infty$}
$$
by the assumption on $D^2u$ in (i). 
The weak maximum principle then implies 
$$
\epsilon\, |x|^{1+\delta} +\max_{\partial B_1}\{-\Psi\} = w(x) \geq -\Psi(x)
\,\,\,\,\mbox{for any $x\in{\mathbb R}^n\backslash B_1$.}
$$
Letting $\epsilon\to 0$, we obtain
$$
\max_{\partial B_1}\{-\Psi\}\geq -\Psi (x)\,\,\,\,\mbox{for any $x\in{\mathbb R}^n\backslash B_1.$}
$$
So $\Psi$ attains its global minimum on ${\mathbb R}^n$ in the closure of $B_1$. Hence $\Psi$ is a constant by applying the strong maximum principle to equation (\ref{psi}).  Now as in the proof of Theorem \ref{maintheorem1},  we  conclude that $u$ must be the quadratic polynomial in the desired form by differentiating equation (\ref{non-euclid}).  Part (i) of Theorem \ref{thm2} is proved. 

If we assume in addition that $u$ is radially symmetric, then $\Psi$ is also radially symmetric and depends only on $|x|$. It follows that $\Psi$ must then attain either a local maximum or a local minimum over any open ball $B$ in ${\mathbb R}^n$. The strong maximum principle then implies $\Psi$ is constant in $B$, and hence in ${\mathbb R}^n$. As before, we conclude that $u$ is quadratic. Part (ii) of Theorem \ref{thm2} is proved. 

\section{Proof of Theorem \ref{thm3}}
Equation \eqref{eqn:10} is in fact closely related to the K\"ahler Ricci flow equation 
\begin{equation}\label{krf}
\frac{\partial g_{i\bar j}}{\partial t} = - R_{i\bar j}.
\end{equation}
Indeed, if $v$ is a strictly-plurisubharmonic solution to \eqref{eqn:10}, then it follows that $u(x,t) = -t\, v\left(\frac{x}{\sqrt{-t}}\right)$ solves the parabolic 
complex Monge-Amp\`ere equation
\begin{equation}
\dfrac{\partial u}{\partial t}= \ln\det(\partial\bar{\partial} u)
\end{equation}
and the K\"ahler metrics $g_{i\bar j}= \partial_i\partial_{\bar j}u$ will evolve according to \eqref{krf}.  Although $g_{i\bar j}$ in general is not a gradient shrinking K\"ahler Ricci soliton, Peng Lu pointed out that it is a shrinking K\"ahler Ricci soliton. 

Now let $v(x)$ and $u(x, t)$ be as above.  Then $v$ satisfies (\ref{eqn:10}), that is, 
\begin{equation}\label{eqn:v}
\ln\partial\bar{\partial}v\left(\frac{x}{\sqrt{-t}}\right)=\frac{1}{2}\frac{x}{\sqrt{-t}}\cdot D v\left(  \frac{x}{\sqrt{-t}
}\right)  - v\left(  \frac{x}{\sqrt{-t}}\right) .
\end{equation}
Introduce the notations
\begin{eqnarray*}
\Phi(x) &=& \ln\det(\partial\bar{\partial} v(x)),\\
g^{i\bar{j}}(x)&=& v^{i\bar{j}}(x). 
\end{eqnarray*}
In the following, we verify 
\begin{eqnarray}
\label{eqn:3}&&\sum_{i,j}g^{i\bar{j}}\Phi_{i\bar{j}}(x) -\frac{1}{2} x\cdot D\Phi(x)=0.
\end{eqnarray} 
 We first calculate 
$$
\Phi_{i} = g^{k\bar l} v_{k\bar{l}i}\,\,\,\,\,\mbox{and}\,\,\,\,\, \Phi_{\bar j} = g^{k\bar l}v_{k\bar{l}\bar{j}}
$$
therefore
$$
x\cdot D\Phi = \sum_k \left( z_k \partial_{z_k}\Phi+\bar{z}_k\partial_{\bar{z}_k}\Phi\right)=\sum_{i,j,k}g^{i\bar j}\left(z_k v_{k\bar{l}i}+\bar{z}_k v_{k\bar{l}\bar{j}}\right).
$$
On the other hand, from (\ref{eqn:v}) we have
\begin{eqnarray*}
\sum_{i,j}g^{i\bar j}\Phi_{i\bar j}&=&\sum_{i,j}g^{i\bar j}\left (\frac{1}{2} x\cdot Dv -v\right)_{i\bar j} \\
&=&\sum_{i,j}g^{i\bar j}\left( \frac{1}{2}\sum_k \left(z_k \partial_{z_k} v +\bar{z}_k \partial_{\bar{z}_k}v   \right)-v \right)_{i\bar j} \\
&=&\frac{1}{2} \sum_{i,j,k}g^{i\bar j} \left(z_k v_{ki\bar j}+\bar{z}_k v_{\bar{k}i\bar j}\right).
\end{eqnarray*}
from which we conclude that equation (\ref{eqn:3}) holds.  

Now take the radial barrier function 
$$
w(r) = \epsilon\, r^{1+\delta} + \max_{\partial B_1} \{ -\Phi\}
$$ 
as in the previous section.  Then we have 
$$
\frac{|x|^2}{2(2m-1+\delta)}\,\Delta w(|x|) -\frac{1}{2}x \cdot Dw(|x|) =0.
$$
Moreover, the assumption on the complex Hessian $\partial\bar{\partial} v$ implies
$$
g^{-1} \leq \frac{2r^2}{2m-1+\delta} I.
$$
Here we assume $\displaystyle{ \partial\bar{\partial}v\geq \frac{(2m-1+\delta)}{2|x|^2}}$ for $|x|>1$ instead of $|x|$ being greater than a large number as in the assumption for simplicity. 
Otherwise, we just replace 1 by the large number, and our arguments go through as well. Note that $\partial\bar{\partial}v>0,$ it 
follows that
\[
\sum_{i,j}g^{i\bar{j}}w_{i\bar{j}}-\frac{1}{2}x\cdot Dw\leq\frac{2\left\vert 
x\right\vert ^{2}}{\left(  2m-1+\delta\right)  }\,\frac{1}{4}\,\Delta
w-\frac{1}{2}x\cdot Dw=0.
\]
Then as in the proof of Theorem \ref{thm2}, the weak maximum principle implies the smooth function $-\Phi$ achieves its global maxima in 
the closure of $B_1$ and the strong maximum principle asserts $\Phi$ is constant. In turn, we conclude that   $v$ is the quadratic polynomial in the claimed form 
by differentiating equation (\ref{eqn:10}) once. Part (i) of Theorem \ref{thm3} is proved. 

The radially symmetric case Part (ii) follows exactly as in the proof of Theorem \ref{thm2}. 

{\it Acknowledgement.}  We are grateful to Peng Lu for the discussion about K\"ahler Ricci solitons.

\bibliographystyle{amsplain}

\end{document}